# SOME APPLICATIONS OF LINEAR ALGEBRA AND GEOMETRY IN REAL LIFE

**Vittoria Bonanzinga**

Mediterranean University of Reggio Calabria, Italy
*vittoria.bonanzinga@unirc.it*

**Abstract**

In this paper, some real-world motivated examples are provided illustrating the power of linear algebra tools as the product of matrices, determinants, eigenvalues and eigenvectors. In this sense, some practical applications related to computer graphics, geometry, areas, volumes are presented, along with some problems connected to sports and investments.

## 1     INTRODUCTION

Numerous studies highlight the need to motivate students to study mathematics by presenting problems related to real life and concrete problems, [9], [10]. The topic described in this work was inspired by the seminars of Alberto Conte and Marina Marchisio held at the University of Messina, May 2017, and at the DIIES Department of the University of Reggio Calabria on the topics of innovation in mathematics teaching, from the conference *Didamatica 2019* (held at the DIGIES Department of the University of Reggio Calabria), the Summer School *Higher Education* (held in Reggio Calabria, June 2019) and the *Aplimat 2020* Conference (held in Bratislava from 2 - 4 February 2020). In literature, there are numerous studies concerning innovative teaching practices in engineering degree courses, Graham, in 2018, [4] presents a substantial 170-page research commissioned and financially supported by the Massachusetts Institute of Technology (MIT), a project on the state of global art in engineering education involving students, faculty, university managers, industrial partners, academic leaders, education researchers, teaching and learning professionals, and national government representatives from around the world who have shared experiences, knowledge and skills. Griese [5] emphasizes the value of students' commitment and motivation for achieving study success.

The tremendous changes in engineering training over the past decade have also affected mathematics training in engineering courses. Furthermore, the development of technology and digital resources has led to new possibilities for engineering courses in which mathematically complex problems are solved with the aid of computers and visualizations and simulations play a central role, [3], [7]. These developments have also changed conditions for teaching and learning. New technologies make it possible to model and solve particularly complex engineering problems with increasingly realistic virtual simulations. The active learning of students on tasks related to real life and self-regulation on individualized learning paths characterize the innovative practices introduced in various universities.

The practical applications of mathematics before being applied to the courses of Geometry for Engineering were successfully tested by the author in the course of Mathematics Fundamentals for the basic training of the Master's Degree Course in Primary Education in the first semester of the academic year 2019/2020, as reported in [2]. To stimulate and support students in the learning process in order to reduce the dropout rate and ensure the achievement of educational success, good didactic planning is essential in which motivation plays a key role. Indeed, students often encounter difficulties in learning mathematics because it is decontextualized, abstract and incomprehensible. The contextualization of real problems that can be solved with the use of linear algebra was a basic element to increase interest and motivation. The educational model presented was inspired by some results obtained within the European Rules-Math Project and was shared by the research group of the University of Turin DELTA, (Digital education for Learning and Teaching Advances).

## 2   APPLICATION IN REAL LIFE OF SOME PROBLEMS OF LINEAR ALGEBRA

Among the innovative elements of the didactic methodologies used there is the contextualization in real life of linear algebra and geometry problems, [9], [10]. Some examples are presented.

### 2.1   Application of matrices in sport

Three friends, Steven, Marc and George decide to train to participate in a triathlon competition (running, swimming and cycling). Every day Steven trains 30 minutes for running, 20 minutes for swimming and 100 minutes for cycling, while Marc trains 25 minutes for running, 30 minutes for swimming and 60 minutes for cycling, and George trains 20 minutes for running, 45 minutes for swimming and 55 minutes for cycling.
  i. Write matrix A, which has the training minutes for each sport of Steven, Marc and George as rows.
  ii. Considering the matrix shown below
$$M = \begin{pmatrix} 10.1 & 9.2 & 12.2 \\ 7.2 & 6.5 & 8.7 \\ 5.3 & 4.6 & 6.4 \end{pmatrix}$$
which represents the calories burned by Steven of 70 kg, Marc of 65 kg, and George of 85 kg in relation to their weight for each minute of training; in the first row there are the calories burned for each minute of training in running, in the second row the calories burned for swimming and in the third for cycling. The first column refers to Steven, the second to Marc and the third to George. Calculate the product rows by columns $A \cdot M$. What do the elements of the main diagonal of the resulting matrix correspond to?
  iii. Which of the three guys, based on their training, consumes the most calories in a day?
  iv. What does the element of the second row and third column of the matrix obtained at point ii) correspond to?

*Solution*
  i. The matrix which represents the training minutes of Steven, Marc and George is the following:
$$A = \begin{pmatrix} 30 & 20 & 100 \\ 25 & 30 & 60 \\ 20 & 45 & 55 \end{pmatrix}.$$
  ii. The product

$$A \cdot M = \begin{pmatrix} 30 & 20 & 100 \\ 25 & 30 & 60 \\ 20 & 45 & 55 \end{pmatrix} \cdot \begin{pmatrix} 10.1 & 9.2 & 12.2 \\ 7.2 & 6.5 & 8.7 \\ 5.3 & 4.6 & 6.4 \end{pmatrix} = C = \begin{pmatrix} 977 & 866 & 1180 \\ 786.5 & 701 & 905 \\ 817.5 & 729.5 & 987.5 \end{pmatrix}.$$

The elements on the main diagonal of the obtained matrix C correspond to the total calories burned for each boy:

$c_{11}$=977 is Steven's total calories burned

$c_{22}$=701 is Marc's total calories burned

$c_{33}$=987.5 is George's total calories burned.

iii. George burns the most calories in one day.

iv. The element of the second row and third column corresponds to the total calories that George would burn if he trained for Marc's minutes as element $c_{23}$ is obtained by multiplying the second row of matrix A containing Marc's training minutes by the third column of the matrix M containing the calories burned by George, in relation to his weight, for each minute of training.

## 2.2 Determinants, areas and volumes

We give a geometric interpretation of the determinant of a matrix of order 2 and of the determinant of a matrix of order 3 and we present some applications to compute areas and volumes using determinants. The absolute value of the determinant of a matrix of order 2 is equal to the area of the parallelogram spanned by the row vectors of the matrix, and the absolute value of the determinant of a matrix of order 3 is simply equal to the volume of the parallelepiped spanned by the row vectors of that matrix, [1], [8]. In dimension 2, we have the following

**Theorem 1** *If A is a matrix of order 2, its rows determine a parallelogram P in $R^2$. The area of the parallelogram P is the absolute value of the determinant of the matrix whose rows are the vectors forming two adjacent sides of the parallelogram:*

$$\text{Area } P = \left| \det \begin{pmatrix} a & b \\ c & d \end{pmatrix} \right| = |ad - bc|.$$

For instance, the area of the parallelogram formed by 2 two-dimensional vectors (1,5) and (6,2) is the following:

$$\text{Area } P = \left| \det \begin{pmatrix} 1 & 5 \\ 6 & 2 \end{pmatrix} \right| = |2 - 30| = 28$$

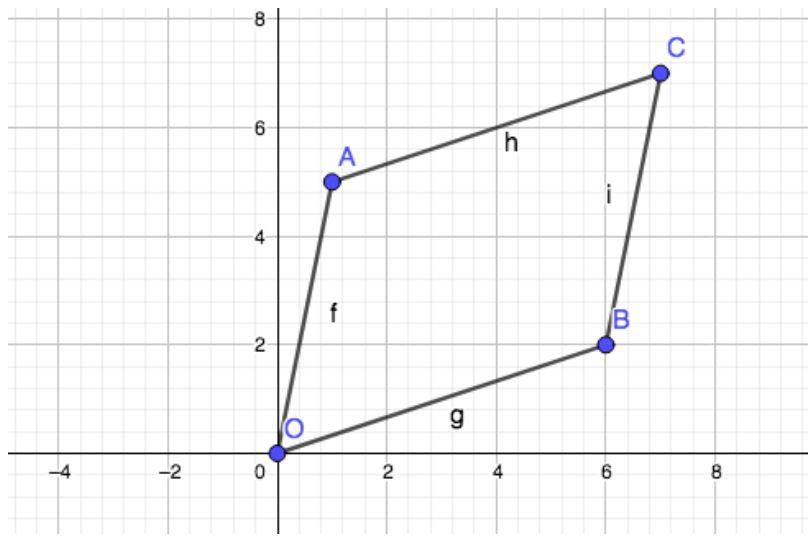

**Figure 1**

squared units.

In general, we consider a parallelogram P and we choose any three vertices of the parallelogram $A(x_A, y_A)$, $B(x_B, y_B)$ and $C(x_C, y_C)$, and two adjacent vectors $CA(x_A-x_C, y_A-y_C)$ and $CB(x_B-x_C, y_B-y_C)$, then

$$Area\ P = \left|\det\begin{pmatrix} x_A - x_C & y_A - y_C \\ x_B - x_C & y_B - y_C \end{pmatrix}\right| = |(x_A - x_C) \cdot (y_B - y_C) - (x_B - x_C) \cdot (y_A - y_C)|.$$

*Area P* can also be written in this elegant way using the coordinates of any three vertices of a parallelogram:

$$Area\ P = \left|\det\begin{pmatrix} x_A & y_A & 1 \\ x_B & y_B & 1 \\ x_C & y_C & 1 \end{pmatrix}\right|.$$

For instance, if A=(1,5), C=(7,7) and B=(6,2) are three vertices of a parallelogram P then

$$Area\ P = \left|\det\begin{pmatrix} 1-7 & 5-7 \\ 6-7 & 2-7 \end{pmatrix}\right| = 28.$$

If $A' = \begin{pmatrix} 4 & 6 \\ 6 & 2 \end{pmatrix}$ then the row vectors of the matrix *A'* generate a parallelogram *P'*, then

$$Area\ P' = \left|\det\begin{pmatrix} 4 & 6 \\ 6 & 2 \end{pmatrix}\right| = |8 - 36| = 28.$$

If we compare the first example, where the row vectors are $R_1(1,5)$ and $R_2(6,2)$ and the second example, where the row vectors are $R'_1(4,6)$ and $R'_2(6,2)$, we can observe that

$$R_1 - R'_1 = R''_1 = (-3,-1)$$

and the corresponding determinant of the matrix $\begin{pmatrix} -3 & -1 \\ 6 & 2 \end{pmatrix}$ is zero. This last result is a consequence of a general property of the determinant of a matrix $A = \begin{pmatrix} a & b \\ c & d \end{pmatrix}$:

$$\det A = \det\begin{pmatrix} R_1 \\ R_2 \end{pmatrix} = \det\begin{pmatrix} R_1 + hR_2 \\ R_2 \end{pmatrix} = \det\begin{pmatrix} R_1 \\ R_2 \end{pmatrix} + \det\begin{pmatrix} hR_2 \\ R_2 \end{pmatrix} \quad \forall h \in \mathbb{R}.$$

From the geometric point of view, the previous result means that the area of the parallelogram does not change if we change the row vectors using elementary transformations. In particular, if we fix the base of the parallelogram and we change a row vector using elementary transformations, the height of the parallelogram does not change and so the area of the parallelogram P and P' coincide. It is possible to see the previous result using GeoGebra, [6], drawing the two parallelograms, the first one P of vertices O(0,0), A(1,5), B(6,2), C(7,7) and the second P' of vertices O(0,0), A'(4,6), B(6,2), C'(10,8):

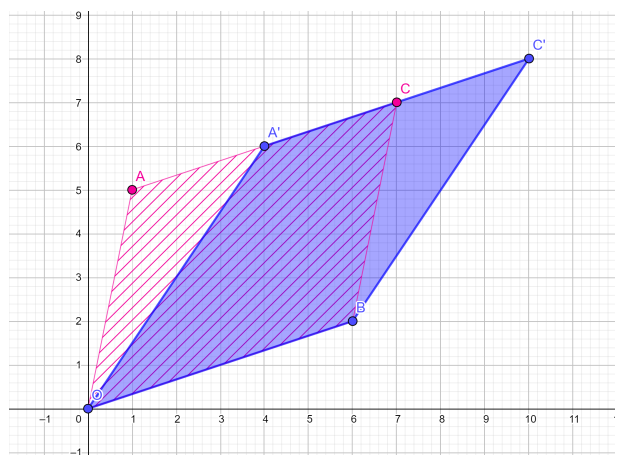

**Figure 2**

By viewing the parallelogram formed by two congruent triangles, we obtain the formula to compute the area of a triangle, if we know the coordinates of its vertices. Precisely, if $A(x_A, y_A)$, $B(x_B, y_B)$ and $C(x_C, y_C)$ are the vertices of a triangle T then

$$\text{Area(T)} = \frac{1}{2} \left| \det \begin{pmatrix} x_A & y_A & 1 \\ x_B & y_B & 1 \\ x_C & y_C & 1 \end{pmatrix} \right|.$$

For instance, if T is a triangle with vertices $A(1,5)$, $B(6,2)$ and $C(7,7)$ then

$$\text{Area(T)} = \frac{1}{2} \left| \det \begin{pmatrix} 1 & 5 & 1 \\ 6 & 2 & 1 \\ 7 & 7 & 1 \end{pmatrix} \right| = 14.$$

In dimension 3, we have the following

**Theorem 2** *If A is a matrix of order 3, its rows, linear independent, determine a parallelepiped P in $R^3$. The volume of the parallelepiped P is the absolute value of the determinant of the matrix A whose rows are three vectors forming three edges of the parallelepiped.*

$$\text{Volume P} = \left| \det \begin{pmatrix} a & b & c \\ a' & b' & c' \\ a'' & b'' & c'' \end{pmatrix} \right|.$$

For instance, the volume of the parallelepiped P determined by the vectors of (1,50), (6,2,0) and (3,2,4) is

$$\text{Volume P} = \left| \det \begin{pmatrix} 1 & 5 & 0 \\ 6 & 2 & 0 \\ 3 & 2 & 4 \end{pmatrix} \right| = 112 \text{ cubic units.}$$

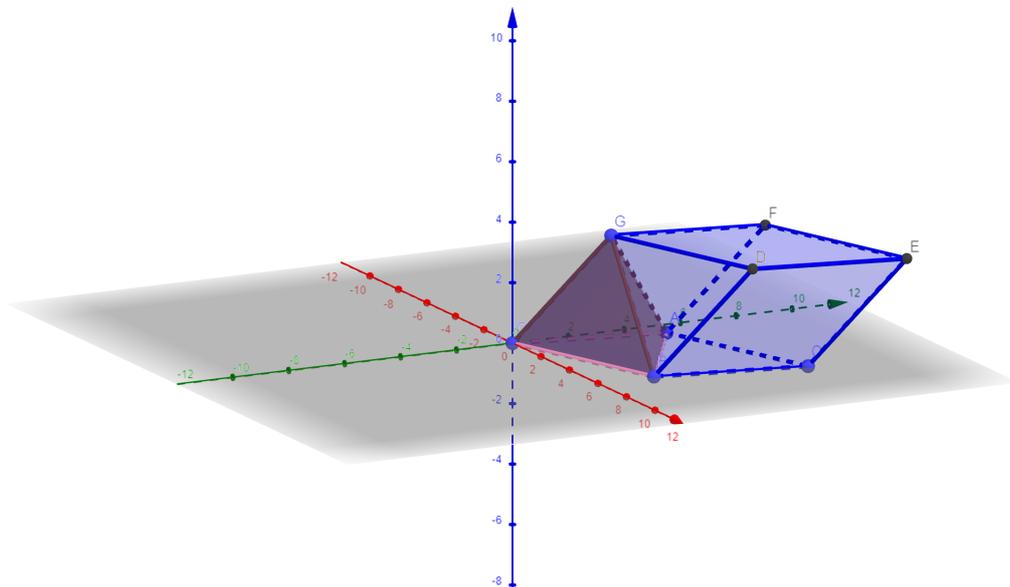

**Figure 3**

As a consequence of the previous result, we can compute the volume of a tetrahedron T, given its vertices, using a determinant. Since the volume of a tetrahedron is $\frac{1}{6}$ of the volume of the corresponding parallelepiped, if $A(x_A, y_A, z_A)$, $B(x_B, y_B, z_B)$, $C(x_C, y_C, z_C)$ and $D(x_D, y_D, z_D)$ are the vertices of a tetrahedron T, then three non-coplanar vectors:

$DA(x_A - x_D, y_A - y_D, z_A - z_D)$,
$DB(x_B - x_D, y_B - y_D, z_B - z_D)$,

$DC(x_C - x_D, y_C - y_D, z_C - z_D)$,
give the desired volume:

$$\text{Vol(T)} = \frac{1}{6}\left|\det\begin{pmatrix} x_A - x_D & y_A - y_D & z_A - z_D \\ x_B - x_D & y_B - y_D & z_B - z_D \\ x_C - x_D & y_C - y_D & z_C - z_D \end{pmatrix}\right|.$$

For instance, if A(1,5,0), B(6,2,0), C(3,2,4), D(0,0,0), are the vertices of a tetrahedron ABCD then the volume is

$$\text{Vol(T)} = \frac{1}{6}\left|\det\begin{pmatrix} 1 & 5 & 0 \\ 6 & 2 & 0 \\ 3 & 2 & 4 \end{pmatrix}\right| = \frac{1}{6}|4(2-30)| = \frac{56}{3}.$$

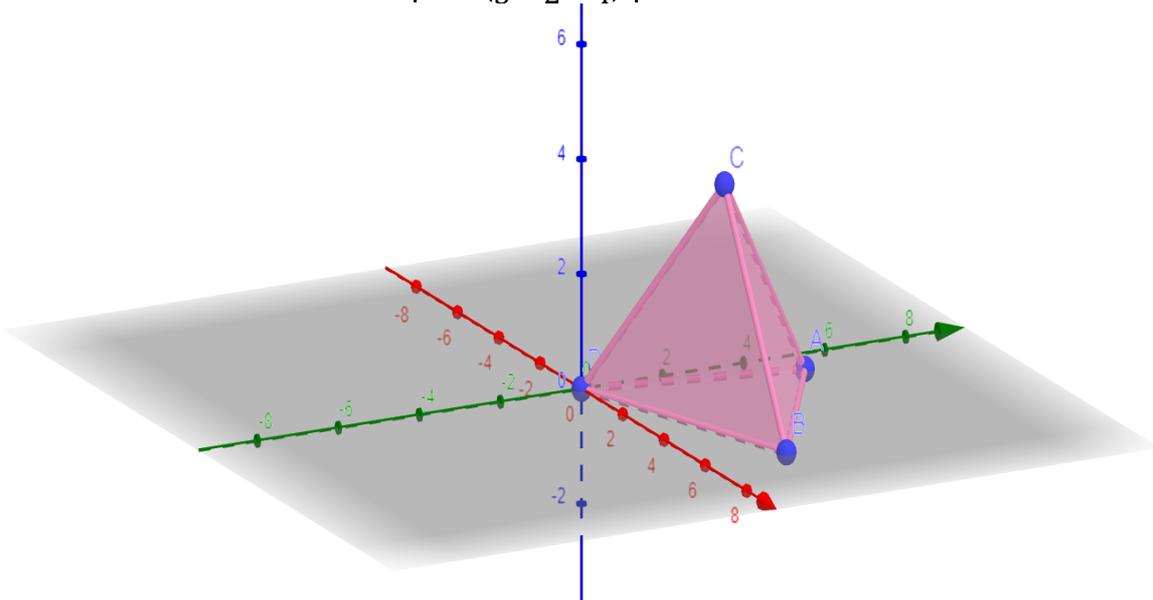

**Figure 4**

**Remark** We can observe in dimension 3 that the volume of a parallelepiped is the area of its base times its height: here the "base" is the parallelogram, determined by row vectors $R_1$ and $R_2$ of matrix A of order 3 and the "height" is the perpendicular distance of $R_3$ from the base. We consider a row replacement of the form $R_3 = R_3 + cR_i$, for i=1,2 and for all $c \in \mathbb{R}$. Translating $R_3$ by a multiple of $R_i$ moves $R_3$ in a direction parallel to the base. This changes neither the base nor the height. Thus, the volume of a parallelepiped is unchanged by row replacements. For instance, if

$$A = \begin{pmatrix} 1 & 5 & 0 \\ 6 & 2 & 0 \\ 3 & 2 & 4 \end{pmatrix}$$

then the corresponding volume of a parallelepiped P spanned by the row vectors of A is given by the absolute value of the determinant of A which is 112. Replacing $R_3$ with $R_3 + 2R_1$ we have

$$A' = \begin{pmatrix} 1 & 5 & 0 \\ 6 & 2 & 0 \\ 5 & 12 & 4 \end{pmatrix}$$

and the corresponding volume of the parallelepiped P' is unchanged. Using GeoGebra, we can visualize parallelepipeds P and P' as follows:

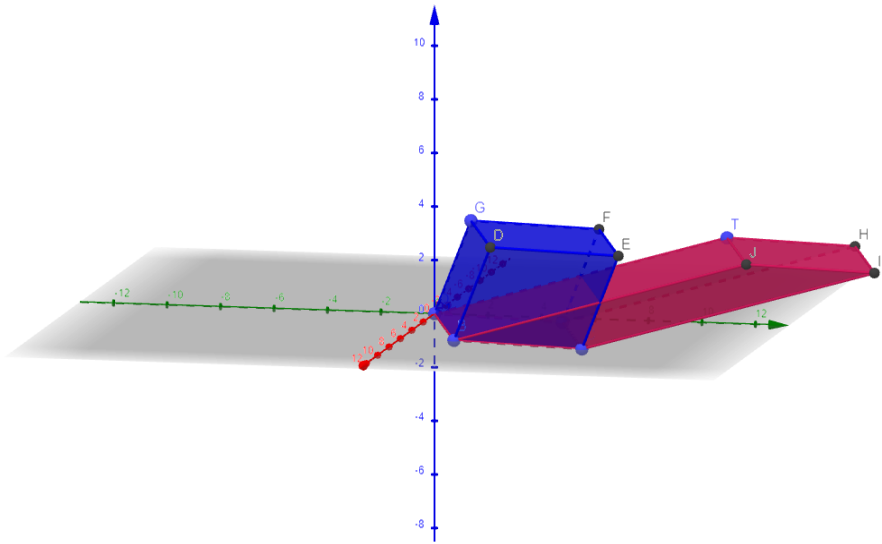

**Figure 5**

### 2.3 Strawberry Fields forever - Investments, eigenvalues and eigenvectors

Let us look at a practical example in the agriculture sphere. A farmer owns land where he grows strawberries. Part of the strawberries are used in three different sectors, part A for the cake production sector, part B for the jam production sector and the last part, part C, is used for the fair in the local village. Section A produces a revenue equal to four times the expenditure made for the strawberries, section B produces revenue of double the expenditure, while part C produces revenue equal to two thirds, since it is based only on offers made during the fair. You want to obtain a revenue proportional to the money invested in strawberry seedlings, whereby, with the revenues, other strawberries are grown which are redistributed equally in the 3 sectors. If the farm initially invests € 4200 in seedlings, how much does it get after a year in the optimal situation?

Let $x$, $y$ and $z$ be the investments in sectors A, B and C respectively.

$$\begin{matrix} A \\ B \\ C \end{matrix} \begin{pmatrix} x \\ y \\ z \end{pmatrix} \xrightarrow{after\ a\ year} \begin{pmatrix} 4x \\ 2y \\ \frac{2z}{3} \end{pmatrix} = \begin{pmatrix} x + 3x \\ y + y \\ \frac{2z}{3} \end{pmatrix} \xrightarrow{redistribution} \begin{pmatrix} x + x + \frac{1}{3}y \\ y + \frac{1}{3}y + x \\ \frac{2z}{3} + x + \frac{1}{3}y \end{pmatrix} = \begin{pmatrix} 2x + \frac{1}{3}y \\ x + \frac{4}{3}y \\ x + \frac{1}{3}y + \frac{2z}{3} \end{pmatrix}.$$

So the transformation of the investment after one year is given by the following linear map:

$$f: \mathbb{R}^3 \to \mathbb{R}^3, with\ f(x,y,z) = \left(2x + \frac{1}{3}y, x + \frac{4}{3}y, x + \frac{1}{3}y + \frac{2}{3}z\right).$$

Let $M = M_{CC}(f) = \begin{pmatrix} 2 & \frac{1}{3} & 0 \\ 1 & \frac{4}{3} & 0 \\ 1 & \frac{1}{3} & \frac{2}{3} \end{pmatrix}$ be the matrix associated with the linear map f with respect to the canonical bases in the domain and in the codomain. Let X be the vector of the initial distribution; we would like the initial distribution after one year to be of the type $\lambda X$, so we look for the eigenvalues and eigenvectors of M. The characteristic polynomial is:

$$det\begin{pmatrix} 2-\lambda & \frac{1}{3} & 0 \\ 1 & \frac{4}{3}-\lambda & 0 \\ 1 & \frac{1}{3} & \frac{2}{3}-\lambda \end{pmatrix} = \left(\frac{2}{3}-\lambda\right)\left(\lambda^2 - \frac{10}{3}\lambda + \frac{7}{3}\right).$$

The eigenspace related to $\lambda_1 = \frac{2}{3}$ is $V_{\lambda_1} = \{(0,0,t): t \in \mathbb{R}\}$.

The eigenspace related to $\lambda_2 = 1$ is $V_{\lambda_2} = \left\{\left(-\frac{1}{3}u, u, 0\right) : u \in \mathbb{R}\right\}$.

The eigenspace related to $\lambda_3 = \frac{7}{3}$ is $V_{\lambda_3} = \left\{\left(r, r, \frac{4}{5}r\right) : r \in \mathbb{R}\right\}$.

Now let us interpret the results in the light of our problem. The eigenvalue $\lambda_1 = \frac{2}{3}$ is a possible solution to our problem: investing all € 4200 in sector C, but this solution would not be a wise choice as every year the investment would decrease more and more. The eigenvalue $\lambda_2 = 1$ does not give any solution to our problem as the eigenvectors $\left(-\frac{1}{3}u, u, 0\right)$ have no interest as we cannot invest a negative amount of money. We consider the eigenvalue $\lambda_3 = \frac{7}{3}$ with the generic eigenvector $\left(r, r, \frac{4}{5}r\right)$. If we invest € 4200, then $r + r + \frac{4}{5}r = €\ 4200$, hence we get $r = €1500$ and $\frac{4}{5}r = €1200$. Thus, the optimal solution is to invest €1500 respectively in sectors A and B and €1200 in sector C; after one year, the amount of money will increase and we will have €3500 in sectors A and B, while €2800 in sector C.

Example 3 provides a simple modelling of a process that can also be applied in the industrial sphere, in the production of goods and in the organization of resources.

### 2.4 Computer graphics and transformations

In the field of computer graphics, the following 2D and 3D transformations are frequently used: translations, rotations, scaling and reflections.

We consider some examples in dimension 2. Given the quadrilateral ABCD obtained by joining the points A(1; 1), B(2; 3), C(4; 3) and D(5; 1), we describe the transformation of the figure ABCD considering the translation of the given figure with a vector of components (3,2), and the rotation of the given figure by an angle $\theta = \frac{\pi}{2}$. We describe these transformations and we represent the figures described above.

We represent the quadrilateral ABCD on an xy plane:

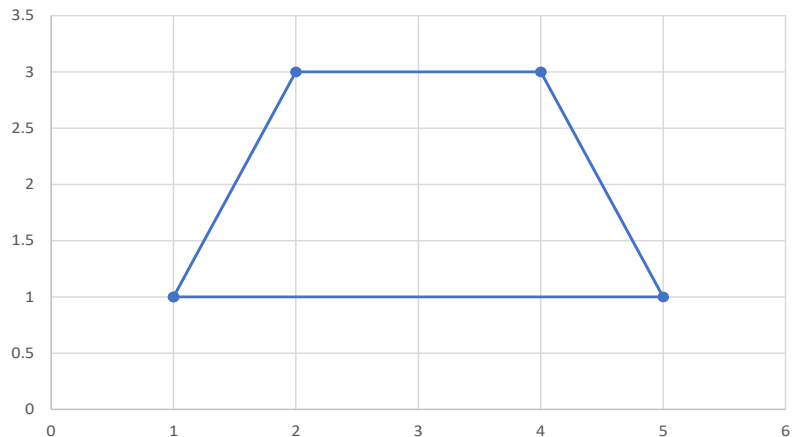

**Figure 6 - Quadrilateral ABCD**

The vector components of this quadrilateral are therefore v₁=(1,1), v₂=(2,3), v₃=(4,3) and v₄=(5,1). A translation is an isometry, that is, a geometric transformation that leaves the distances unchanged by moving all the points by a fixed distance in the same direction. The translation $T_V$ is the application $T_V: \mathbb{R}^2 \to \mathbb{R}^2$ and therefore:

$$T_v\begin{pmatrix}x\\y\end{pmatrix} = \begin{pmatrix}x+p\\y+q\end{pmatrix}.$$

where p and q are the components of the translation vector. If p=3 and q=2 then by applying the translation to the vectors v₁=(1,1), v₂=(2,3), v₃=(4,3) and v₄=(5,1) we obtain the vectors of components: (4,3), (5,5), (7,5) e (8,3) respectively, illustrated in Figure 7

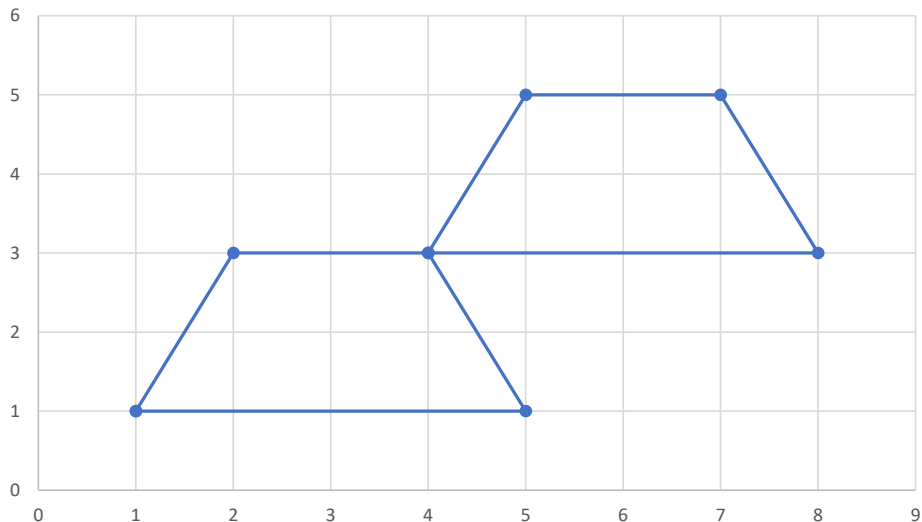

**Figure 7 – Quadrilateral ABCD and its translation**

Rotation in the plane is a map of $\mathbb{R}^2 \to \mathbb{R}^2$ described by the following function:

$$R(\theta) = \begin{pmatrix}\cos\theta & -\sin\theta\\ \sin\theta & \cos\theta\end{pmatrix}\begin{pmatrix}x\\y\end{pmatrix}$$

If $\theta = \frac{\pi}{2}$, then applying the rotation to the vectors v₁=(1,1), v₂=(2,3), v₃=(4,3) and v₄=(5,1) we obtain respectively the vectors of components: (-1,1), (-3,2), (-3,4) and (-1,5) described in the following figure:

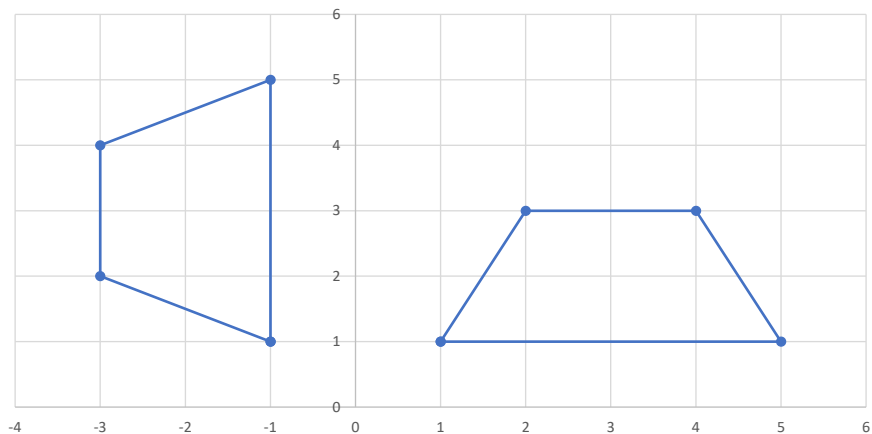

**Figure 8 – Quadrilateral ABCD and its rotation**

We consider some examples in dimension 3. Given the tetrahedron ABCD obtained by joining the points A(1,5,0), B(6,2,0), C(3,2,4), D(0,0,0), we describe the transformation of the figure ABCD considering the translation of the given figure with a vector of components (4,3,-2), the rotation of the given figure by an angle $\theta = \frac{\pi}{2}$ about the z-axis and the reflection of the tetrahedron ABCD through the xy-plane. We describe these transformations and we represent the figures described above. The tetrahedron ABCD is represented in Figure 4.

The vector components of the tetrahedron ABCD are therefore $v_1=(1,5,0)$, $v_2=(6,2,0)$, $v_3=(3,2,4)$ and $v_4=(0,0,0)$. In dimension 3, the translation $T_V$ is the application $T_V: \mathbb{R}^3 \rightarrow \mathbb{R}^3$ and therefore:

$$T_v \begin{pmatrix} x \\ y \\ z \end{pmatrix} = \begin{pmatrix} x+p \\ y+q \\ z+r \end{pmatrix}$$

where p, q and r are the components of the translation vector. If p=4, q=3 and r=-2 then by applying the translation to the vectors $v_1=(1,5,0)$, $v_2=(6,2,0)$, $v_3=(3,2,4)$ and $v_4=(0,0,0)$ we obtain the vectors of components: (5,8,-2), (10,5,-2), (7,5,2) and (4,3,-2) respectively, illustrated in Figure 9.

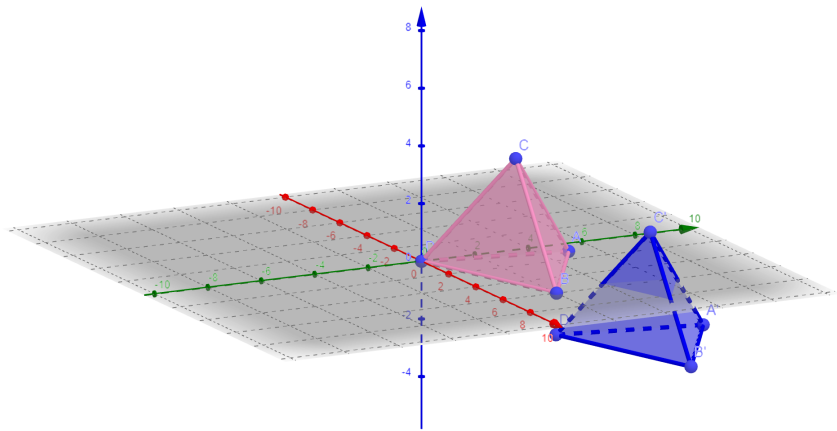

**Figure 9- Tetrahedron ABCD and its translation**

Rotation in the space of the given figure by an angle $\theta$ about the z-axis is a map of $\mathbb{R}^3 \rightarrow \mathbb{R}^3$ described by the following function:

$$R(\theta) = \begin{pmatrix} \cos\theta & -\sin\theta & 0 \\ \sin\theta & \cos\theta & 0 \\ 0 & 0 & 1 \end{pmatrix} \begin{pmatrix} x \\ y \\ z \end{pmatrix}.$$

If $\theta = \frac{\pi}{2}$, $R\left(\frac{\pi}{2}\right) = \begin{pmatrix} 0 & -1 & 0 \\ 1 & 0 & 0 \\ 0 & 0 & 1 \end{pmatrix} \begin{pmatrix} x \\ y \\ z \end{pmatrix} = \begin{pmatrix} -y \\ x \\ z \end{pmatrix}$ then applying the rotation to the vectors $v_1=(1,5,0)$, $v_2=(6,2,0)$, $v_3=(3,2,4)$ and $v_4=(0,0,0)$ we obtain respectively the vectors of components: (-5,1,0), (-2,6,0), (-2,3,4) and (0,0,0) described in the following figure:

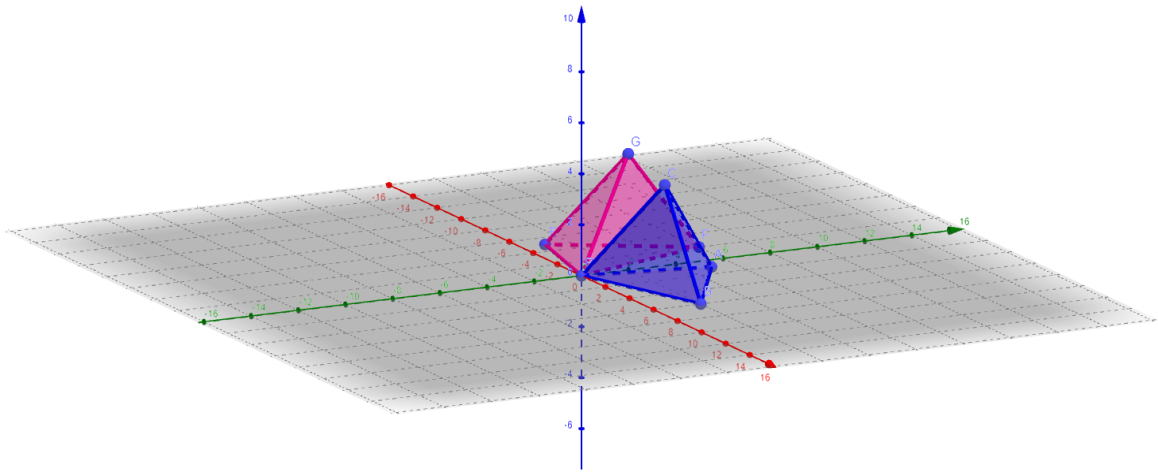

**Figure 10 – Tetrahedron ABCD and its rotation**

The reflection of the tetrahedron ABCD through the xy-plane is a map of $\mathbb{R}^3 \rightarrow \mathbb{R}^3$ described by the following function: $s(x, y, z) = (x, y, -z)$ then applying the reflection to the vectors $v_1$=(1,5,0), $v_2$=(6,2,0), $v_3$=(3,2,4) and $v_4$=(0,0,0) we obtain respectively the vectors of components: (1,5,0), (6,2,0), (3,2,-4), and (0,0,0) described in the following figure:

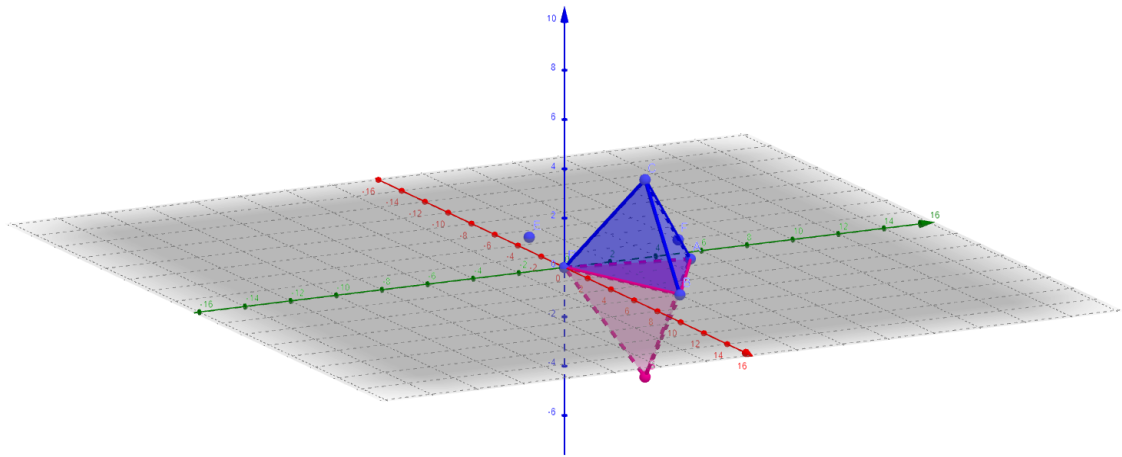

**Figure 11 - Tetrahedron ABCD and its reflection**

**Acknowledgement**

Vittoria Bonanzinga is member of the National Group for Algebraic and Geometric Structures and their applications (GNSAGA/INdAM).